\documentclass{amsart}

\input{prepictex}
\input{pictex}
\input{postpictex}

\begin{document}

\def\sbullet{\raise2.3pt\hbox{$_\bullet$}}
\def\bddots{\kern-10pt\lower7pt\hbox{.}\kern-7pt\lower11pt\h
box{.}
\kern-7pt\lower15pt\hbox{.}}
\def\uptoLambda(#1){\Lambda^{\le #1}}

\def\CC{\mathbb C}
\def\NN{\mathbb N}
\def\RR{\mathbb R}
\def\TT{\mathbb T}
\def\ZZ{\mathbb Z}

\def\Bb{{\mathcal B}}
\def\Ff{{\mathcal F}}
\def\Hh{{\mathcal H}}
\def\Kk{{\mathcal K}}
\def\Oo{{\mathcal O}}
\def\Ss{{\mathcal S}}

\def\Ad{\operatorname{Ad}}
\def\Aut{\operatorname{Aut}}
\def\clsp{\overline{\lsp}}
\def\End{\operatorname{End}}
\def\id{\operatorname{id}}
\def\Isom{\operatorname{Isom}}
\def\Ker{\operatorname{Ker}}
\def\lsp{\operatorname{span}}
\def\Hom{\operatorname{Hom}}
\def\Obj{\operatorname{Obj}}

\theoremstyle{plain}
\newtheorem{theorem}{Theorem}[section]
\newtheorem*{theorem*}{Theorem}
\newtheorem*{prop*}{Proposition}
\newtheorem{cor}[theorem]{Corollary}
\newtheorem{lemma}[theorem]{Lemma}
\newtheorem{prop}[theorem]{Proposition}
\theoremstyle{remark}
\newtheorem{rmk}[theorem]{Remark}
\newtheorem{rmks}[theorem]{Remarks}
\newtheorem*{aside}{Aside}
\newtheorem*{note}{Note}
\newtheorem{comment}[theorem]{Comment}
\newtheorem{example}[theorem]{Example}
\newtheorem{examples}[theorem]{Examples}
\theoremstyle{definition}
\newtheorem{dfn}[theorem]{Definition}
\newtheorem{dfns}[theorem]{Definitions}
\newtheorem{notation}[theorem]{Notation}

\numberwithin{equation}{section}

\title
{Higher-rank graphs and their $C^*$-algebras}
\author{Iain Raeburn}
\author{Aidan Sims}
\author{Trent Yeend}
\address{Department of Mathematics  \\
       University of Newcastle\\  NSW  2308\\ AUSTRALIA}
\date{21 June 2001}
\thanks{This research was supported by the Australian
Research Council.}

\begin{abstract}
We consider the higher-rank graphs introduced by Kumjian and
Pask as models for higher-rank Cuntz-Krieger algebras. We
describe a variant of the Cuntz-Krieger relations which
applies to graphs with sources, and describe a local
convexity condition which characterises the higher-rank
graphs that admit a nontrivial Cuntz-Krieger family. We then
prove versions of the uniqueness theorems and
classifications of ideals for the $C^*$-algebras generated
by Cuntz-Krieger families.
\end{abstract}

\maketitle

\section{Introduction}

The $C^*$-algebras of higher-rank graphs of Kumjian and Pask
\cite{KP} generalise the higher-rank Cuntz-Krieger algebras
of Robertson and Steger \cite{RS1, RS2, RS3} in the same way
that the $C^*$-algebras of infinite graphs generalise the
original Cuntz-Krieger algebras \cite{C, CK}. In \cite{KP},
Kumjian and Pask analysed higher-rank graph algebras using a
groupoid model like that used in \cite{KPRR} and \cite{KPR}
to analyse graph algebras. The results in \cite{KPRR} and
\cite{KPR} were sharpened in \cite{BPRS} using a direct
analysis based on the  original arguments used by Cuntz and
Krieger in \cite{C} and \cite{CK}; the analysis of
\cite{BPRS} applies to the algebras of quite general
row-finite graphs, and in particular the graphs can have
sinks or sources.

Here we carry out a direct analysis of the $C^*$-algebras of
row-finite higher-rank graphs. One interesting new feature
is the difficulty in extending results to higher-rank graphs
with sources: the paths $\lambda$ in higher-rank graphs have
degrees $d(\lambda)$ in $\NN^k$ rather than lengths in
$\NN,$ and vertices may receive edges of some degrees and
not of others. To overcome this difficulty we modify the
Cuntz-Krieger relation to ensure that Cuntz-Krieger algebras
have a spanning family of the usual sort, and identify a
local convexity condition under which a higher-rank graph
admits a nontrivial Cuntz-Krieger family. We then prove
versions of the gauge-invariant uniqueness theorem and the
Cuntz-Krieger uniqueness theorem for locally convex
higher-rank graphs, extending the results of \cite{KP}, and
use them to investigate the ideal structure.

The Cuntz-Krieger relation of \cite{KP} involves the spaces
$\Lambda^p$ of paths of degree $p\in \NN^k$. Our key
technical innovation is the introduction of path spaces
$\uptoLambda(p)$ consisting of the paths $\lambda$ with
$d(\lambda)\leq p$ which cannot be extended to paths
$\lambda\mu$ with $d(\lambda\mu)\leq p$; the key
Lemmas~\ref{lambda inclusion} and \ref{orth range proj} say
that the spaces $\uptoLambda(p)$ have combinatorial
properties like those of the spaces $\Lambda^p$, and ensure
that the $C^*$-algebras behave like Cuntz-Krieger algebras
(see Proposition~\ref{theorem 4.11}). These new path spaces
would also have simplified the analysis of the core in the
$C^*$-algebras of graphs with sinks in \cite[\S2]{BPRS}.
Indeed, the $\uptoLambda(q)$ notation works so smoothly that
arguments sometimes appear deceptively easy.

\section{Higher-rank graphs}

\begin{dfns} Given $k\in\NN$, a {\it graph of rank $k$} (or
$k$-graph) $(\Lambda,d)$ consists of a countable category
$\Lambda=(\Obj(\Lambda),\Hom(\Lambda),r,s)$ together with a
functor $d:\Lambda\to\NN^k$, called the {\it degree map},
which satisfies the {\it factorisation property}: for every
$\lambda\in\Lambda$ and $m,n\in\NN^k$ with $d(\lambda)=m+n$,
there are unique elements $\mu,\nu\in\Lambda$ such that
$\lambda =\mu\nu$, $d(\mu)=m$ and $d(\nu)=n$. Elements
$\lambda\in\Lambda$ are called {\it paths}. For $m\in\NN^k$
and $v\in\Obj(\Lambda)$, we define
$\Lambda^m:=\{\lambda\in\Lambda: d(\lambda)=m\}$ and
$\Lambda^m(v):=\{\lambda\in\Lambda^m: r(\lambda)=v\}$. A
morphism between two $k$-graphs $(\Lambda_1,d_1)$ and
$(\Lambda_2,d_2)$ is a functor $f:\Lambda_1\to\Lambda_2$
which respects the degree maps. $(\Lambda,d)$ is {\it
row-finite} if for each $v\in\Obj(\Lambda)$ and $m\in\NN^k$,
the set $\Lambda^m(v)$ is finite; $(\Lambda,d)$ has {\it no
sources} if $\Lambda^m(v)\neq\emptyset$ for all
$v\in\Obj(\Lambda)$ and $m\in\NN^k$. \end{dfns}

The factorisation property says that there is a unique path
of degree $0$ at each vertex, and hence allows us to
identify $\Obj(\Lambda)$ with $\Lambda^0$.

\begin{examples} \label{k-graphs}
(i) Let $E$ be a directed graph, and $l:E^*\to \NN$ the
length function on the path space. Then $(E^*,l)$ is a
$1$-graph.

(ii) Let $k \in\NN$, let $m\in(\NN\cup\{\infty\})^k$, and
define a partial ordering on $\NN^k$ by $m\le n \iff m_i\le
n_i$ for all $i$. $(\Omega_{k,m},d)$ is the $k$-graph with
category $\Omega_{k,m}$ defined by
\begin{gather*}
\Obj(\Omega_{k,m}) := \{p\in \NN^k : p\le m\}, \\
\Hom(\Omega_{k,m}) := \{(p,q) \in \Obj(\Omega_{k,m}) \times
\Obj(\Omega_{k,m}) : p \le q\},
\end{gather*}
$r(p,q):=p$, $s(p,q):=q$, and degree map $d(p,q):=q-p$.

When each $m_i=\infty$, the resulting $k$-graph is the main
example $\Omega_k$ used in \cite{KP}.

When $m\in\NN^k$, the resulting $k$-graphs are important
because every path $\lambda$ of degree $m$ in a $k$-graph
$\Lambda$ determines a graph morphism
$x_\lambda:\Omega_{k,m}\to \Lambda$: set
$x_\lambda(p,q):=\lambda''$ where
$\lambda=\lambda'\lambda''\lambda'''$ and $d(\lambda')=p$,
$d(\lambda''')=m-q$. Indeed, this sets up a bijection
between $\Lambda^m$ and the graph morphisms
$x:\Omega_{k,m}\to\Lambda$.
\end{examples}

To visualise a $k$-graph, we draw its \emph{$1$-skeleton},
which is the graph with vertex set $\Lambda^0$, edge set
$\bigcup^k_{i=1} \Lambda^{e_i}$, range and source maps
inherited from $\Lambda$, and with the edges of different
degrees distinguished using $k$ different colours. (In the
pictures here, we imagine that dashed lines are red and
dotted lines are blue.) For example, the $1$-skeleton of
$\Omega_{2,(3,2)}$ is

\begin{equation}\label{omega(3,2)}
\beginpicture
\setcoordinatesystem units <1cm,1cm>
\put{$\bullet$} at 1 -1
\put{$\bullet$} at 2 -1
\put{$\bullet$} at 3 -1
\put{$\bullet$} at 4 -1
\put{$\bullet$} at 1 0
\put{$\bullet$} at 2 0
\put{$\bullet$} at 3 0
\put{$\bullet$} at 4 0
\put{$\bullet$} at 1 1
\put{$\bullet$} at 2 1
\put{$\bullet$} at 3 1
\put{$\bullet$} at 4 1
\setlinear
\plot 1.1 -1 1.9 -1 /
\plot 2.1 -1 2.9 -1 /
\plot 3.1 -1 3.9 -1 /
\plot 1.1 0 1.9 0 /
\plot 2.1 0 2.9 0 /
\plot 3.1 0 3.9 0 /
\plot 1.1 1 1.9 1 /
\plot 2.1 1 2.9 1 /
\plot 3.1 1 3.9 1 /
\arrow <0.15cm> [0.25,0.75] from  1.52 -1 to 1.48 -1
\arrow <0.15cm> [0.25,0.75] from  2.52 -1 to 2.48 -1
\arrow <0.15cm> [0.25,0.75] from  3.52 -1 to 3.48 -1
\arrow <0.15cm> [0.25,0.75] from  1.52 0 to 1.48 0
\arrow <0.15cm> [0.25,0.75] from  2.52 0 to 2.48 0
\arrow <0.15cm> [0.25,0.75] from  3.52 0 to 3.48 0
\arrow <0.15cm> [0.25,0.75] from  1.52 1 to 1.48 1
\arrow <0.15cm> [0.25,0.75] from  2.52 1 to 2.48 1
\arrow <0.15cm> [0.25,0.75] from  3.52 1 to 3.48 1
\setdashes
\plot 1 -.9 1 -.1 /
\plot 2 -.9 2 -.1 /
\plot 3 -.9 3 -.1 /
\plot 1 .1 1 .9 /
\plot 2 .1 2 .9 /
\plot 3 .1 3 .9 /
\plot 4 -.9 4 -.1 /
\plot 4 .1 4 .9 /
\arrow <0.15cm> [0.25,0.75] from  1 -.48 to 1 -.52
\arrow <0.15cm> [0.25,0.75] from  1 .52 to 1 .48
\arrow <0.15cm> [0.25,0.75] from  2 -.48 to 2 -.52
\arrow <0.15cm> [0.25,0.75] from  2 .52 to 2 .48
\arrow <0.15cm> [0.25,0.75] from  3 -.48 to 3 -.52
\arrow <0.15cm> [0.25,0.75] from  3 .52 to 3 .48
\arrow <0.15cm> [0.25,0.75] from  4 -.48 to 4 -.52
\arrow <0.15cm> [0.25,0.75] from  4 .52 to 4 .48
\put{$(0,0)$} at .5 -1.1
\put{$(3,2)$} at 4.5 1.1
\put{$g$} at 4.2 .5
\put{$e$} at 3.45 -.15
\put{$f$} at 2.8 .5
\put{$h$} at 3.45 1.2
\put{$v$} at 2.82 -.18
\endpicture
\end{equation}

\noindent
Because the edges represent morphisms in a category, we
write $eg$ for the path in the $1$-skeleton which consists
of $g$ followed by $e$.

The $1$-skeleton of a $k$-graph does not always suffice to
determine the $k$-graph: we have to say how the edges in
$\Lambda^{e_i}$ fit together to give elements of
$\Lambda^n$. We interpret elements of $\Lambda^n$ as
commuting diagrams of shape $n$ in which the morphisms
correspond to edges in the given $1$-skeleton. Thus, for
example, in a $2$-graph with $1$-skeleton

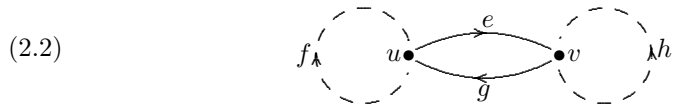
\begin{equation}\label{niceex}
\beginpicture
\setcoordinatesystem units <1cm,1cm>
\put{$\bullet$} at 4 0
\put{$\bullet$} at 6 0
\setdashes
\circulararc 320 degrees from 4 .2 center at 3.4 0
\arrow <0.15cm> [0.25,0.75] from 2.77  .01 to 2.77 .03
\circulararc 320 degrees from 6 -.2 center at 6.6 0
\arrow <0.15cm> [0.25,0.75] from 7.227 .075 to 7.226 .095
\setsolid
\setquadratic
\plot 4.1 .08 5 .3 5.9 .08 /
\arrow <0.15cm> [0.25,0.75] from 4.98 .3 to 5.02 .3
\setquadratic
\plot 4.1 -.07 5 -.3 5.9 -.07 /
\arrow <0.15cm> [0.25,0.75] from 4.92 -.3 to 4.88 -.3
\put{$u$} at 3.8  0
\put{$v$} at 6.2 0
\put{$f$} at 2.6 0
\put{$e$} at 5.05 .45
\put{$g$} at 5 -.5
\put{$h$} at 7.4 .1
\endpicture
\end{equation}

\noindent
where the dashed  lines have degree $(0,1)$, the unique
example of a $(3,1)$ path $\lambda$ with $r(\lambda)=u$ and
$s(\lambda)=v$ is

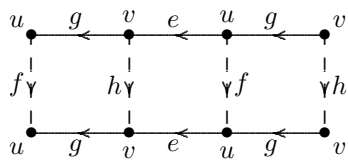
\[
\beginpicture
\setcoordinatesystem units <1.3cm,1.3cm>
\put{$\bullet$} at 1 0
\put{$\bullet$} at 2 0
\put{$\bullet$} at 3 0
\put{$\bullet$} at 4 0
\put{$\bullet$} at 1 1
\put{$\bullet$} at 2 1
\put{$\bullet$} at 3 1
\put{$\bullet$} at 4 1
\setlinear
\setsolid
\plot 1.05 0 1.95 0 /
\plot 2.05 0 2.95 0 /
\plot 1.05 1 1.95 1 /
\plot 2.05 1 2.95 1 /
\plot 3.05 0 3.95 0 /
\plot 3.05 1 3.95 1 /
\arrow <0.15cm> [0.25,0.75] from  1.52 0 to 1.48 0
\arrow <0.15cm> [0.25,0.75] from  2.52 0 to 2.48 0
\arrow <0.15cm> [0.25,0.75] from  1.52 1 to 1.48 1
\arrow <0.15cm> [0.25,0.75] from  2.52 1 to 2.48 1
\arrow <0.15cm> [0.25,0.75] from  3.52 0 to 3.48 0
\arrow <0.15cm> [0.25,0.75] from  3.52 1 to 3.48 1
\setdashes
\plot 1 .1 1 1 /
\plot 2 .1 2 1 /
\plot 3 .1 3 1 /
\plot 4 .1 4 1 /
\arrow <0.15cm> [0.25,0.75] from  1 .45 to 1 .41
\arrow <0.15cm> [0.25,0.75] from  2 .45 to 2 .41
\arrow <0.15cm> [0.25,0.75] from  3 .45 to 3 .41
\arrow <0.15cm> [0.25,0.75] from  4 .45 to 4 .41
\put{$u$} at .85 -.15
\put{$u$} at .85 1.15
\put{$u$} at 3 -.2
\put{$u$} at 3 1.2
\put{$v$} at 2 1.2
\put{$v$} at 2 -.2
\put{$v$} at 4.15 1.15
\put{$v$} at 4.15 -.15
\put{$g$} at 1.45 -.15
\put{$g$} at 1.45 1.15
\put{$e$} at 2.45 -.12
\put{$e$} at 2.45 1.15
\put{$f$} at .85 .5
\put{$h$} at 1.85 .5
\put{$f$} at 3.15 .5
\put{$g$} at 3.45 -.15
\put{$g$} at 3.45 1.15
\put{$h$} at 4.15 .5
\endpicture
\]

\noindent
 From such a picture we can read off the factorisations of
$\lambda$: $\lambda=gegh=gefg=gheg=fgeg$.

When $k=2$, it suffices to specify the factorisations of
paths $ef$ of length $2$ in the $1$-skeleton for which $e$
and $f$ have different colours. Any collection $S$ of
squares which contains each such bi-coloured path exactly
once determines a unique $2$-graph $\Lambda$ with the given
$1$-skeleton and $\Lambda^{(1,1)}=S$ (see \cite[\S6]{KP});
there may be no such collection, or there may be many. For
the $1$-skeleton in (\ref{niceex}), the factorisation
property implies that $\Lambda^{(1,1)}$ consists of the two
squares
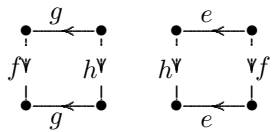
\[
\beginpicture
\setcoordinatesystem units <1cm,1cm>
\put{$\bullet$} at 1 0
\put{$\bullet$} at 2 0
\put{$\bullet$} at 1 1
\put{$\bullet$} at 2 1
\setlinear
\setsolid
\plot 1.1 0 1.9 0 /
\plot 1.1 1 1.9 1 /
\arrow <0.15cm> [0.25,0.75] from  1.52 0 to 1.48 0
\arrow <0.15cm> [0.25,0.75] from  1.52 1 to 1.48 1
\setdashes
\plot 1 .1 1 .9 /
\plot 2 .1 2 .9 /
\arrow <0.15cm> [0.25,0.75] from  1 .52 to 1 .48
\arrow <0.15cm> [0.25,0.75] from  2 .52 to 2 .48
\put{$g$} at 1.4 -.2
\put{$g$} at 1.4 1.2
\put{$f$} at .85 .5
\put{$h$} at 1.85 .5
\put{$\bullet$} at 3 0
\put{$\bullet$} at 4 0
\put{$\bullet$} at 3 1
\put{$\bullet$} at 4 1
\setlinear
\setsolid
\plot 3.1 0 3.9 0 /
\plot 3.1 1 3.9 1 /
\arrow <0.15cm> [0.25,0.75] from  3.52 0 to 3.48 0
\arrow <0.15cm> [0.25,0.75] from  3.52 1 to 3.48 1
\setdashes
\plot 3 .1 3 .9 /
\plot 4 .1 4 .9 /
\arrow <0.15cm> [0.25,0.75] from  3 .52 to 3 .48
\arrow <0.15cm> [0.25,0.75] from  4 .52 to 4 .48
\put{$e$} at 3.4 -.2
\put{$e$} at 3.4 1.2
\put{$h$} at 2.85 .5
\put{$f$} at 4.15 .5
\endpicture
\]

\noindent
and hence there is exactly one $2$-graph with this
$1$-skeleton. However, if we add one extra edge to the
$1$-skeleton in (\ref{niceex}), we have to make a choice.
For example, in the $1$-skeleton
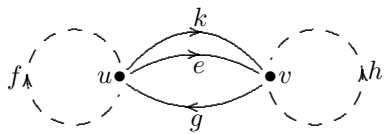
\[
\beginpicture
\setcoordinatesystem units <1cm,1cm>
\put{$\bullet$} at 4 0
\put{$\bullet$} at 6 0
\setdashes
\circulararc 320 degrees from 4 .2 center at 3.4 0
\arrow <0.15cm> [0.25,0.75] from 2.77  .01 to 2.77 .03
\circulararc 320 degrees from 6 -.2 center at 6.6 0
\arrow <0.15cm> [0.25,0.75] from 7.227 .075 to 7.226 .095
\setsolid
\setquadratic
\plot 4.15 .05 5 .3 5.85 .05 /
\arrow <0.15cm> [0.25,0.75] from 5.02 .3 to 5.06 .3
\plot 4.1 .1 5 .6 5.9 .1 /
\arrow <0.15cm> [0.25,0.75] from 5.04 .6 to 5.08 .6
\setquadratic
\plot 4.1 -.1 5 -.4 5.9 -.1 /
\arrow <0.15cm> [0.25,0.75] from 4.92 -.4 to 4.88 -.4
\put{$u$} at 3.8  0
\put{$v$} at 6.2 0
\put{$f$} at 2.6 0
\put{$e$} at 5.05 .13
\put{$k$} at 5.07 .8
\put{$g$} at 5.02 -.6
\put{$h$} at 7.4 .1
\endpicture
\]
there are four possible bi-coloured paths from $u$ to $v$,
and we have to decide how to pair these off into paths of
degree $(1,1)$: either
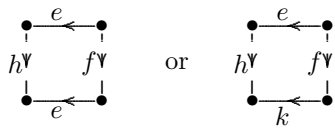
\[
\beginpicture
\setcoordinatesystem units <1cm,1cm>
\put{$\bullet$} at 1 0
\put{$\bullet$} at 2 0
\put{$\bullet$} at 1 1
\put{$\bullet$} at 2 1
\setlinear
\setsolid
\plot 1.1 0 1.9 0 /
\plot 1.1 1 1.9 1 /
\arrow <0.15cm> [0.25,0.75] from  1.52 0 to 1.48 0
\arrow <0.15cm> [0.25,0.75] from  1.52 1 to 1.48 1
\setdashes
\plot 1 .1 1 .9 /
\plot 2 .1 2 .9 /
\arrow <0.15cm> [0.25,0.75] from  1 .52 to 1 .48
\arrow <0.15cm> [0.25,0.75] from  2 .52 to 2 .48
\put{$e$} at 1.4 -.15
\put{$e$} at 1.4 1.15
\put{$h$} at .85 .5
\put{$f$} at 1.83 .5
\put{or} at 3 .5
\put{$\bullet$} at 4 0
\put{$\bullet$} at 5 0
\put{$\bullet$} at 4 1
\put{$\bullet$} at 5 1
\setlinear
\setsolid
\plot 4.1 0 4.9 0 /
\plot 4.1 1 4.9 1 /
\arrow <0.15cm> [0.25,0.75] from  4.52 0 to 4.48 0
\arrow <0.15cm> [0.25,0.75] from  4.52 1 to 4.48 1
\setdashes
\plot 4 .1 4 .9 /
\plot 5 .1 5 .9 /
\arrow <0.15cm> [0.25,0.75] from  4 .52 to 4 .48
\arrow <0.15cm> [0.25,0.75] from  5 .52 to 5 .48
\put{$k$} at 4.4 -.2
\put{$e$} at 4.4 1.15
\put{$h$} at 3.85 .5
\put{$f$} at 4.83 .5
\endpicture
\]
is a path of degree $(1,1)$, and once we have decided which,
the other pairing is determined by the factorisation
property.

For $k>2$, a collection $S$ of squares may not be the set of
paths of degree $(1,1)$ for any $k$-graph with  the given
$1$-skeleton. However, \cite[Theorem~2.1]{FS} tells us that
it suffices to know that for every tri-coloured path $efg$
in the $1$-skeleton, the six squares on the sides of the
cube
\bigskip
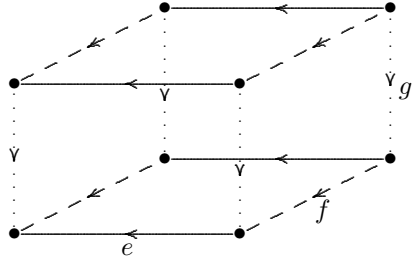
\[
\beginpicture
\setcoordinatesystem units <1cm,1cm>
\put{$\bullet$} at 1 -1
\put{$\bullet$} at 4 -1
\put{$\bullet$} at 6  0
\put{$\bullet$} at 3 0
\put{$\bullet$} at 1 1
\put{$\bullet$} at 4 1
\put{$\bullet$} at 6 2
\put{$\bullet$} at 3 2
\setlinear
\plot 1.1 -1 3.9 -1 /
\plot 3.1 0 5.9 0 /
\plot 1.1 1 3.9 1 /
\plot 3.1 2 5.9 2 /
\arrow <0.15cm> [0.25,0.75] from 2.52 -1 to 2.48 -1
\arrow <0.15cm> [0.25,0.75] from 2.52 1 to 2.48 1
\arrow <0.15cm> [0.25,0.75] from 4.52 0 to 4.48 0
\arrow <0.15cm> [0.25,0.75] from 4.52 2 to 4.48 2
\setdashes
\plot 4.1 -.95 5.9 -.05 /
\plot 1.1 -.95 2.9 -.05 /
\plot 4.1 1.05 5.9 1.95 /
\plot 1.1 1.05 2.9 1.95 /
\arrow <0.15cm> [0.25,0.75] from 5.02 -.49 to 4.98 -.51
\arrow <0.15cm> [0.25,0.75] from 2.02 -.49 to 1.98 -.51
\arrow <0.15cm> [0.25,0.75] from 2.02 1.51 to 1.98 1.49
\arrow <0.15cm> [0.25,0.75] from 5.02 1.51 to 4.98 1.49
\setdots
\plot 1 -.9 1 .9 /
\plot 4 -.9 4 .9 /
\plot 3 .1 3 1.9 /
\plot 6 .1 6 1.9 /
\setsolid
\arrow <0.15cm> [0.25,0.75] from 1 .02 to 1 -.02
\arrow <0.15cm> [0.25,0.75] from 4 -.18 to 4  -.22
\arrow <0.15cm> [0.25,0.75] from 3 .82 to 3 .78
\arrow <0.15cm> [0.25,0.75] from 6 1.02 to 6 .98
\put{$e$} at 2.5 -1.2
\put{$f$} at 5.1 -.7
\put{$g$} at 6.2 .9
\endpicture
\]
\bigskip

\noindent
give a well-defined path of degree $(1,1,1)$. More
precisely, we need to know that the path $g^2f^2e^2$ with
reverse colouring obtained by successively filling in the
three visible squares agrees with the path $g_2f_2e_2$
obtained by filling in the three invisible squares:
\bigskip
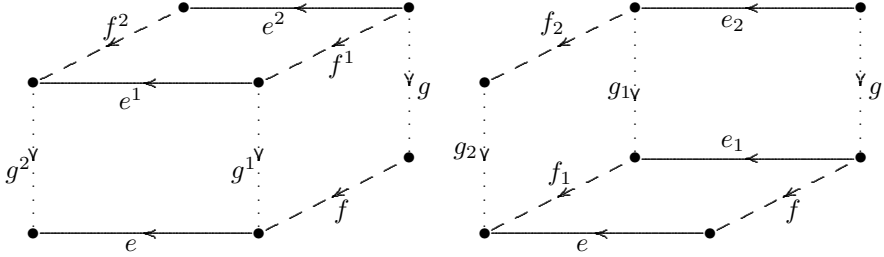
\[
\beginpicture
\setcoordinatesystem units <1cm,1cm>
\put{$\bullet$} at 1 -1
\put{$\bullet$} at 4 -1
\put{$\bullet$} at 6  0
\put{$\bullet$} at 1 1
\put{$\bullet$} at 4 1
\put{$\bullet$} at 6 2
\put{$\bullet$} at 3 2
\setlinear
\plot 1.1 -1 3.9 -1 /
\plot 1.1 1 3.9 1 /
\plot 3.1 2 5.9 2 /
\arrow <0.15cm> [0.25,0.75] from 2.52 -1 to 2.48 -1
\arrow <0.15cm> [0.25,0.75] from 2.52 1 to 2.48 1
\arrow <0.15cm> [0.25,0.75] from 4.52 2 to 4.48 2
\setdashes
\plot 4.1 -.95 5.9 -.05 /
\plot 4.1 1.05 5.9 1.95 /
\plot 1.1 1.05 2.9 1.95 /
\arrow <0.15cm> [0.25,0.75] from 5.02 -.49 to 4.98 -.51
\arrow <0.15cm> [0.25,0.75] from 2.02 1.51 to 1.98 1.49
\arrow <0.15cm> [0.25,0.75] from 5.02 1.51 to 4.98 1.49
\setdots
\plot 1 -.9 1 .9 /
\plot 4 -.9 4 .9 /
\plot 6 .1 6 1.9 /
\setsolid
\arrow <0.15cm> [0.25,0.75] from 1 .02 to 1 -.02
\arrow <0.15cm> [0.25,0.75] from 4  .02 to 4  -.02
\arrow <0.15cm> [0.25,0.75] from 6 1.02 to 6 .98
\put{$e$} at 2.3 -1.15
\put{$f$} at 5.1 -.7
\put{$g$} at 6.2 .9
\put{$g^1$} at 3.8 -.15
\put{$f^1$} at 5.1 1.25
\put{$g^2$} at .8 -.15
\put{$e^1$} at 2.3 .8
\put{$f^2$} at 2.1 1.75
\put{$e^2$} at 4.2 1.8
\put{$\bullet$} at 7 -1
\put{$\bullet$} at 10 -1
\put{$\bullet$} at 12  0
\put{$\bullet$} at 9 0
\put{$\bullet$} at 7 1
\put{$\bullet$} at 12 2
\put{$\bullet$} at 9 2
\setsolid
\setlinear
\plot 7.1 -1 9.9 -1 /
\plot 9.1 0 11.9 0 /
\plot 9.1 2 11.9 2 /
\arrow <0.15cm> [0.25,0.75] from 8.52 -1 to 8.48 -1
\arrow <0.15cm> [0.25,0.75] from 10.52 0 to 10.48 0
\arrow <0.15cm> [0.25,0.75] from 10.52 2 to 10.48 2
\setdashes
\plot 10.1 -.95 11.9 -.05 /
\plot 7.1 -.95 8.9 -.05 /
\plot 7.1 1.05 8.9 1.95 /
\arrow <0.15cm> [0.25,0.75] from 11.02 -.49 to 10.98 -.51
\arrow <0.15cm> [0.25,0.75] from 8.02 -.49 to 7.98 -.51
\arrow <0.15cm> [0.25,0.75] from 8.02 1.51 to 7.98 1.49
\setdots
\plot 7 -.9 7 .9 /
\plot 9 .1 9 1.9 /
\plot 12 .1 12 1.9 /
\setsolid
\arrow <0.15cm> [0.25,0.75] from 7 .02 to 7 -.02
\arrow <0.15cm> [0.25,0.75] from 9 .82 to 9 .78
\arrow <0.15cm> [0.25,0.75] from 12 1.02 to 12 .98
\put{$e$} at 8.3 -1.15
\put{$f$} at 11.1 -.7
\put{$g$} at 12.2 .9
\put{$f_1$} at 8 -.2
\put{$e_1$} at 10.3 .2
\put{$g_1$} at 8.8 .9
\put{$g_2$} at 6.75 .1
\put{$f_2$} at 7.9 1.8
\put{$e_2$} at 10.3 1.8
\endpicture
\]
(In the left-hand diagram, we first use the right hand face
to determine $g^1f^1$, then use the front to find $g^2e^1$,
and then the top to find $f^2e^2$; in the right-hand
diagram, we use the bottom first.) A family of squares which
contains each bi-coloured path exactly once and which
satisfies this condition on cubes determines a unique
$k$-graph (see Theorem~2.1 and Remark~2.3 in \cite{FS}).

\section{Higher-rank graphs and their $C^*$-algebras}

For a row-finite $k$-graph $\Lambda$ with no sources the
authors of \cite{KP} define a Cuntz-Krieger $\Lambda$-family
to be a family of partial isometries
$\{s_\lambda:\lambda\in\Lambda\}$ such that
$\{s_v:v\in\Lambda^0\}$ are mutually orthogonal projections,
$s_\lambda s_\mu=s_{\lambda\mu}$ for all
$\lambda,\mu\in\Lambda$ with $r(\mu)=s(\lambda)$,
$s^*_\lambda s_\lambda=s_{s(\lambda)}$ for all
$\lambda\in\Lambda$, and
\begin{equation} \label{KP def}
s_v=\sum_{\lambda\in\Lambda^m(v)}s_\lambda s^*_\lambda
\quad\text{for all
$v\in\Lambda^0$ and $m\in\NN^k$.}
\end{equation}
When $\Lambda$ has sources there is trouble with
relation~(\ref{KP def}) in the sense that $\Lambda^m(v)$ may
be nonempty for some values of $m$ and empty for others. If
a vertex $v\in\Lambda^0$ receives no edges we could just
impose no relation for that vertex as is done for directed
graphs in \cite{KPR}; if $v$ receives edges of some degrees
but not others, then we have to change (\ref{KP def}) in
more subtle ways. The obvious strategy is to observe that
when $\Lambda^{e_i}(v) \neq \emptyset$ for all
$i\in\{1,\dots,k\}$ and $v \in \Lambda^0$, (\ref{KP def}) is
equivalent to (\ref{revised KP def}), and then to replace
(\ref{KP def}) with
\begin{equation}\label{revised KP def}
s_v=\sum_{\lambda\in\Lambda^{e_i}(v)}s_\lambda s^*_\lambda
\quad\text{for
$v\in\Lambda^0$ and $1\le i\le k$ with
$\Lambda^{e_i}(v)\neq\emptyset$.}
\end{equation}
While this works for a large class of $k$-graphs (see
Proposition~\ref{equiv rel}), in general there are problems.
We consider the following key example.
\begin{equation}\label{locally concave}
\beginpicture
\setcoordinatesystem units <1.3cm,1.3cm>
\put{$\bullet$} at 1 -.5
\put{$\bullet$} at 2 -.5
\put{$\bullet$} at 1 .5
\setlinear
\setsolid
\plot 1 -.5 2 -.5 /
\arrow <0.15cm> [0.25,0.75] from  1.52 -.5 to 1.48 -.5
\setdashes
\plot 1 -.4 1 .5 /
\arrow <0.15cm> [0.25,0.75] from  1 -.05 to 1 -.09
\put{$e$} at 1.45 -.6
\put{$f$} at .85 0
\put{$v$} at .87 -.59
\put{$w$} at 2.15 -.6
\put{$z$} at .9 .62
\endpicture
\end{equation}
In the $2$-graph of (\ref{locally concave}),
relation~\eqref{revised KP def} would say
$s_es^*_e=s_v=s_fs^*_f$. But then $s^*_es_f$ would be a
partial isometry with source projection $s_w$ and range
projection $s_z$, and consequently would not be expressible
as a sum of partial isometries of the form $s_\lambda
s^*_\mu$; thus the $C^*$-algebra generated by $\{s_f, s_e,
s_v\}$ would not look like a Cuntz-Krieger algebra. This
problem does not arise in the $2$-graph given by
(\ref{omega(3,2)}): the compositions $eg$ and $fh$ define
the same path in $\Lambda^{(1,1)}$, so $s_es_g=s_fs_h$, and
$s^*_es_f=s_gs^*_h$.

Our adaptation of (\ref{KP def}) retains relations for each
$m\in\NN^k$, but involves sums over paths which extend as
far as possible in the direction $m$. Formally, we
introduce:

\begin{dfn}\label{uptolambda def}
Let $(\Lambda,d)$ be a $k$-graph. For $q\in\NN^k$ and
$v\in\Lambda^0$ we define
\[
\uptoLambda(q):=\{\lambda\in\Lambda: d(\lambda)\le q, \text{
and }\Lambda^{e_i}(s(\lambda))=\emptyset \text{ when }
d(\lambda)+e_i\le q\},
\]
and
\[
\uptoLambda(q)(v):=\{\lambda\in\uptoLambda(q):r(\lambda)=v\}
\]
\end{dfn}

\begin{rmks}
Notice that $\uptoLambda(q)(v)$ is never empty: if there are
no nontrivial paths of degree less than or equal to $q$,
then $\uptoLambda(q)(v)=\{v\}$; in particular, if
$r^{-1}(v)=\emptyset$, then $\uptoLambda(q)(v)=\{v\}$ for
all $q\in\NN^k$. The sets $\uptoLambda(q)$ and
$\uptoLambda(q)(v)$ are used in arguments where the
$\Lambda^q$ and $\Lambda^q(v)$ may have been used in
\cite{KP}; when $\Lambda$ has no sources, $\uptoLambda(q) =
\Lambda^q$.
\end{rmks}

\begin{dfn} \label{k-graph algebra dfn}
Let $\Lambda$ be a row-finite $k$-graph. A {\it
Cuntz-Krieger $\Lambda$-family} in a $C^*$-algebra $B$
consists of a family of partial isometries $\{ s_\lambda
:\lambda\in\Lambda\}$ satisfying the {\it Cuntz-Krieger
relations}:
\begin{itemize}
\item[(1)]
$\{s_v:v\in\Lambda^0\}$ is a family of mutually orthogonal
projections,
\item[(2)]
$s_{\lambda\mu}=s_\lambda s_\mu$ for all
$\lambda,\mu\in\Lambda$ with $s(\lambda)=r(\mu)$,
\item[(3)]
$s^*_\lambda s_\lambda = s_{s(\lambda)}$,
\item[(4)]
$\displaystyle{s_v=\sum_{\lambda\in\uptoLambda(m)(v)}s_\lambda s^*_\lambda}$ for all
$v\in\Lambda^0$ and $m\in\NN^k$.
\end{itemize}
\end{dfn}

\begin{examples}\label{zero rep ex}
For the $2$-graph given by (\ref{omega(3,2)}), the relations
at $v$ say $s_v=s_es^*_e$ for $m=(1,0)$, $s_v=s_fs^*_f$ for
$m=(0,1)$, and $s_v=s_{eg}s^*_{eg}=s_{fh}s^*_{fh}$ for
$m=(1,1)$. For any $m\in\NN^2$ with $m\ge (1,1)$, the
relation at $v$ for $m$ reduces to that for $(1,1)$.

In the $2$-graph given by (\ref{locally concave}), the
relation $s_v=s_es^*_e+s_fs^*_f$ for $m=(1,1)$ combines with
$s_es^*_e=s_v=s_fs^*_f$ to force everything to be zero.
\end{examples}

The following proposition shows that our Cuntz-Krieger
relations yield the usual type of spanning family (see
Remark~\ref{post theorem 4.11 rmk}~(1)).

\begin{prop} \label{theorem 4.11}
Let $(\Lambda,d)$ be a row-finite $k$-graph and let $\{
s_\lambda:\lambda\in\Lambda\}$ be a Cuntz-Krieger
$\Lambda$-family. Then for $\lambda,\mu\in\Lambda$ and
$q\in\NN^k$ with $d(\lambda),d(\mu)\le q$ we have
\begin{equation} \label{spanning family}
s^*_\lambda s_\mu = \sum_{\substack{\lambda\alpha = \mu\beta
\\ \lambda\alpha\in\uptoLambda(q)}}s_\alpha s^*_\beta.
\end{equation}
\end{prop}

To prove this, we need some properties of $\uptoLambda(q)$.

\begin{lemma} \label{lambda inclusion}
Let $(\Lambda,d)$ be a $k$-graph,
$\lambda\in\uptoLambda(m)$, and
$\alpha\in\uptoLambda(n)(s(\lambda))$. Then
$\lambda\alpha\in\uptoLambda(m+n)$.
\end{lemma}

\begin{proof}
We know $d(\lambda\alpha)\le m+n$. Suppose there exists $i$
such that $d(\lambda\alpha)+e_i\le m+n$. If
$d(\alpha)+e_i\le n$, then
$\Lambda^{e_i}(s(\lambda\alpha))=\Lambda^{e_i}(s(\alpha))=\emptyset$, so suppose not. Then
$\langle d(\alpha),e_i\rangle=\langle n,e_i\rangle$, so
$d(\lambda\alpha)+e_i\le m+n$ implies that
$d(\lambda)+e_i\le m$. But $\lambda\in\uptoLambda(m)$, so
$\Lambda^{e_i}(s(\lambda))=\emptyset$. Now suppose that
there exists
$\beta\in\Lambda^{e_i}(s(\lambda\alpha))=\Lambda^{e_i}(s(\alpha))$. Then by the factorisation
property there exist
$\mu_1,\mu_2\in\Lambda$ such that $\mu_1\mu_2=\alpha\beta$
and $d(\mu_1)=e_i$. But then
$\mu_1\in\Lambda^{e_i}(s(\lambda))$, a contradiction.
Therefore $\Lambda^{e_i}(s(\lambda\alpha))=\emptyset$, and
$\lambda\alpha\in\uptoLambda({m+n})$.
\end{proof}

\begin{lemma} \label{orth range proj}
Let $(\Lambda,d)$ be a row-finite $k$-graph and let $\{
s_\lambda:\lambda\in\Lambda\}$ be a Cuntz-Krieger
$\Lambda$-family. Then for $q\in\NN^k$ and
$\lambda,\mu\in\uptoLambda(q)$, $s^*_\lambda
s_\mu=\delta_{\lambda,\mu}s_{s(\lambda)}$.
\end{lemma}

\begin{proof}
The Cuntz-Krieger relations (1) and (4) tell us that the
projections $\{ s_\alpha
s^*_\alpha:\alpha\in\uptoLambda(q)\}$ are mutually
orthogonal. Cuntz-Krieger relations (2) and (3) then give
\begin{flalign*}
&&s^*_\lambda s_\mu =(s^*_\lambda s_\lambda)s^*_\lambda
s_\mu(s^*_\mu s_\mu)
&=s^*_\lambda(s_\lambda s^*_\lambda)(s_\mu s^*_\mu)s_\mu
=\delta_{\lambda,\mu}s_{s(\lambda)}. &\qed
\end{flalign*}
\renewcommand{\qed}{}\end{proof}

\begin{proof}[Proof of Proposition~\ref{theorem 4.11}]
\begin{flalign*}
&&s^*_\lambda s_\mu &= s_{s(\lambda)}s^*_\lambda s_\mu
s_{s(\mu)}
                 \quad\text{using Definition~\ref{k-graph
algebra dfn} (2)}& \\
&&&=
\Bigg(\sum_{\alpha\in\uptoLambda({q-d(\lambda)})(s(\lambda))
}
              s_\alpha s^*_\alpha\Bigg)s^*_\lambda s_\mu
    \Bigg(\sum_{\beta\in\uptoLambda({q-d(\mu)})(s(\mu))}
              s_\beta s^*_\beta\Bigg)& \\
&&&\hskip5cm\text{using Definition~\ref{k-graph algebra dfn}
(4)}& \\
&&&=
\Bigg(\sum_{\alpha\in\uptoLambda({q-d(\lambda)})(s(\lambda))
}
              s_\alpha s^*_{\lambda\alpha}\Bigg)
    \Bigg(\sum_{\beta\in\uptoLambda({q-d(\mu)})(s(\mu))}
              s_{\mu\beta}s^*_\beta\Bigg)& \\
&&&\hskip5cm\text{using Definition~\ref{k-graph algebra dfn}
(2)}& \\
&&&= \sum_{\alpha\in\uptoLambda({q-d(\lambda)})(s(\lambda))}
    \sum_{\beta\in\uptoLambda({q-d(\mu)})(s(\mu))}
        s_\alpha s^*_{\lambda\alpha}s_{\mu\beta}s^*_\beta& \\
&&&= \sum_{\substack{\lambda\alpha = \mu\beta \\
\lambda\alpha\in\uptoLambda(q)(r(\lambda))}}
        s_\alpha s_{s(\alpha)} s^*_\beta
                 \quad\text{ by Lemma~\ref{lambda inclusion}
and Lemma~\ref{orth range proj}}& \\
&&&= \sum_{\substack{\lambda\alpha = \mu\beta \\
\lambda\alpha\in\uptoLambda(q)(r(\lambda))}}
        s_\alpha s^*_\beta
                 \quad\text{using Definition~\ref{k-graph
algebra dfn} (2)} &\qed
\end{flalign*}
\renewcommand{\qed}{}
\end{proof}

\begin{rmks}\label{post theorem 4.11 rmk}
(1) Proposition~\ref{theorem 4.11} implies that
\[
C^*(\{s_\lambda : \lambda \in \Lambda\})=\clsp\{s_\alpha
s^*_\beta: s(\beta)=s(\alpha)\}.
\]

(2) When $q=d(\lambda)\vee d(\mu)$, if
$\lambda\alpha=\mu\beta$ and
$\lambda\alpha\in\uptoLambda(q)(r(\lambda))$, then
$d(\lambda\alpha)=d(\mu\beta)=q$. Notice that if we take
$\lambda=\mu$ in Proposition~\ref{theorem 4.11}, then the
lemma reduces to relation~(4) of Definition~\ref{k-graph
algebra dfn} at the vertex $s(\lambda)$. Hence (4) could be
replaced with relation \eqref{spanning family} from
Proposition~\ref{theorem 4.11}: this relation would be
harder to verify in examples, but might provide more
insight.
\end{rmks}

Given a row-finite $k$-graph $(\Lambda,d)$, there is a
$C^*$-algebra $C^*(\Lambda)$ generated by a universal
Cuntz-Krieger $\Lambda$-family
$\{s_\lambda:\lambda\in\Lambda\}$ (see \cite[\S 1]{B}); in
other words, for each Cuntz-Krieger $\Lambda$-family
$\{t_\lambda:\lambda\in\Lambda\}$, there is a homomorphism
$\pi:C^*(\Lambda)\to C^*(\{t_\lambda\})$ such that
$\pi(s_\lambda)=t_\lambda$ for every $\lambda\in\Lambda$.
Contrary to our experience with directed graphs and their
$C^*$-algebras, it is not straightforward to construct
Cuntz-Krieger $\Lambda$-families $\{t_\lambda\}$ in which
all the partial isometries $t_\lambda$ are nonzero; in fact,
as we saw in Examples~\ref{zero rep ex}, the $2$-graph of
(\ref{locally concave}) admits no such families. We will
show that the existence of nontrivial Cuntz-Krieger families
is characterised by a local convexity condition on the
$k$-graph.

\begin{dfn}
A $k$-graph $(\Lambda,d)$ is {\it locally convex} if, for
all $v \in \Lambda^0,\ i,j \in \{1, \dots, k\}$ with $i
\not= j,\ \lambda \in \Lambda^{e_i}(v)$ and $\mu \in
\Lambda^{e_j}(v)$, $\Lambda^{e_j}(s(\lambda))$ and
$\Lambda^{e_i}(s(\mu))$ are nonempty.
\end{dfn}

\begin{rmk}
The $2$-graph of (\ref{locally concave}) is not locally
convex since $r^{-1}(u)=r^{-1}(w)=\emptyset$. Every
$1$-graph is locally convex, as is every higher-rank graph
without sources.
\end{rmk}

\begin{prop}\label{equiv rel}
Let $(\Lambda,d)$ be a locally convex row-finite $k$-graph.
Then the Cuntz-Krieger relation~{\rm (4)} of
Definition~\ref{k-graph algebra dfn} is equivalent to
\textnormal{(\ref{revised KP def}).}
\end{prop}

The crucial idea in the proof of Proposition~\ref{equiv rel}
is that, when the $k$-graph is locally convex, the
factorisation property of paths extends to elements of
$\uptoLambda(m)(v)$. This is not the case for the $2$-graph
of (\ref{locally concave}): the path $f$ is in
$\uptoLambda({(1,1)})(v)$, but does not factor as
$\lambda'\lambda''$ with
$\lambda'\in\uptoLambda({(1,0)})(v)$.

\begin{lemma}\label{spadesuit lemma}
Let $(\Lambda,d)$ be a $k$-graph, and suppose that
$(\Lambda,d)$ is locally convex. Then for all
$v\in\Lambda^0$, $m\in\NN^k\setminus\{0\}$, and
$j\in\{1,\dots,k\}$ with $\langle m,e_j\rangle\ge 1$, we
have
\[
\uptoLambda(m)(v)=\{\lambda'\lambda''\in\Lambda:\lambda'\in
\uptoLambda({m-e_j})( v) \text{
and }
\lambda''\in\uptoLambda({e_j})(s(\lambda'))\}.
\]
\end{lemma}

\begin{proof}
Fix $j\in\{1,\dots,k\}$ with $\langle m,e_j\rangle\ge 1$;
there is at least one such $j$ since $m\neq 0$. If
$\lambda'\lambda''\in\Lambda$ satisfies
$\lambda'\in\uptoLambda({m-e_j})(v)$ and
$\lambda''\in\uptoLambda({e_j})(s(\lambda'))$, then by
Lemma~\ref{lambda inclusion}
$\lambda'\lambda''\in\uptoLambda(m)(v)$, so we have the
containment
\[
\uptoLambda(m)(v)\supseteq\{\lambda'\lambda''\in\Lambda:\lambda'\in\uptoLambda({m-e_j})(v)
\text{ and }
\lambda''\in\uptoLambda({e_j})(s(\lambda'))\}.
\]

Now suppose $\lambda\in\uptoLambda(m)(v)$. Since
$d(\lambda)\le m$, we must have:
\begin{itemize}
\item[(1)]
$\langle d(\lambda),e_j\rangle < \langle m,e_j\rangle$, or
\item[(2)]
$\langle d(\lambda),e_j\rangle = \langle m,e_j\rangle$.
\end{itemize}
First suppose (1) holds. Then $d(\lambda)\le m-e_j$, and
hence $\lambda\in\uptoLambda({m-e_j})$. Also
$\Lambda^{e_j}(s(\lambda))=\emptyset$ (since
$\lambda\in\uptoLambda(m)(v)$, and hence
$\uptoLambda({e_j})(s(\lambda))=\{s(\lambda)\}$). Taking
$\lambda'=\lambda$ and $\lambda''=s(\lambda)$ gives
$\lambda=\lambda'\lambda''$ with
$\lambda'\in\uptoLambda({m-e_j})(v)$ and
$\lambda''\in\uptoLambda({e_j})(s(\lambda'))$. Now suppose
(2) holds. Then we can factorise $\lambda=\lambda'\lambda''$
with $d(\lambda'')=e_j$, so
$\lambda''\in\uptoLambda({e_j})(s(\lambda'))$. We claim that
$\lambda'\in\uptoLambda({m-e_j})(v)$. To see this, suppose
that $\lambda'\notin\uptoLambda({m-e_j})(v)$. Then there
exists $i$ such that $d(\lambda')+e_i\le m-e_j$ and
$\Lambda^{e_i}(s(\lambda'))\neq\emptyset$, say
$\alpha\in\Lambda^{e_i}(s(\lambda'))$. By (2) we know that
$\langle d(\lambda'),e_j\rangle=\langle m-e_j,e_j\rangle$,
so $i\neq j$. Since $\Lambda$ is locally convex there is a
$\beta\in\Lambda^{e_i}(s(\lambda''))$, but this implies
$d(\lambda\beta)\le m$, a contradiction since
$\lambda\in\uptoLambda(m)(v)$. Hence
$\lambda=\lambda'\lambda''$ with
$\lambda'\in\uptoLambda({m-e_j})(v)$ and
$\lambda''\in\uptoLambda({e_j})(s(\lambda'))$.
\end{proof}

\begin{rmk}
The $k$-graphs $(\Lambda,d)$ studied in \cite{KP} have no
sources; that is, $\lambda^m(v)\neq\emptyset$ for all
$v\in\Lambda^0$ and $m\in\NN^k$. Then $\uptoLambda(m)(v) =
\Lambda^m(v)$, and Lemma~\ref{spadesuit lemma} just says
that a path $\lambda\in\Lambda^m(v)$ can be factorised into
$\lambda=\lambda'\lambda''$ with $d(\lambda')=m-e_j$ and
$d(\lambda'')=e_j$. Thus local convexity ensures that
$\uptoLambda(m)(v)$ has factorisation properties like those
of $\Lambda^m(v)$.
\end{rmk}

\begin{proof}[Proof of Proposition~\ref{equiv rel}]
Property (\ref{revised KP def}) merely consists of specific
cases of Cuntz-Krieger relation (4), so suppose
\eqref{revised KP def} holds. Define a map $l:\NN^k\to \NN$
by $l(m)=\sum^k_{i=1} m_i$. We prove (4) by induction on
$l(n).$

The $n=0$ case is trivial since it amounts to the tautology
$s_v=s_v$, and the $n=1$ case follows directly from
(\ref{revised KP def}) since $l^{-1}(1)=\{ e_1,\dots,e_k\}$.
Suppose (4) holds for all $m$ such that $l(m)\le n$. Let
$m\in \NN^k$ satisfy $l(m)=n+1$ and choose $j$ such that
$m_j\ge1$. Using Lemma~\ref{spadesuit lemma} we have
\begin{align*}
\sum_{\lambda\in\uptoLambda(m)(v)}s_\lambda s^*_\lambda
&= \sum_{\lambda'\in\uptoLambda({m-e_j})(v)}
\sum_{\lambda''\in\uptoLambda({e_j})(s(\lambda'))}
        s_{\lambda'\lambda''}s^*_{\lambda'\lambda''} \\
&= \sum_{\lambda'\in\uptoLambda({m-e_j})(v)}
        s_{\lambda'}
\Bigg(\sum_{\lambda''\in\uptoLambda({e_j})(s(\lambda'))}
s_{\lambda''}s^*_{\lambda''}\Bigg)
        s^*_{\lambda'} \\
&= \sum_{\lambda'\in\uptoLambda({m-e_j})(v)}
        s_{\lambda'}s_{s(\lambda')}s^*_{\lambda'}
         \quad\text{by (\ref{revised KP def})}\\
&= \sum_{\lambda'\in\uptoLambda({m-e_j})(v)}
        s_{\lambda'} s^*_{\lambda'} \\
&= s_v \quad\text{by inductive hypothesis}.
\end{align*}
Hence (4) holds whenever $l(m) = n+1$.
\end{proof}

Kumjian and Pask use the infinite path space
$\Lambda^\infty$ of a $k$-graph $\Lambda$ with no sources to
produce a Cuntz-Krieger $\Lambda$-family of nonzero partial
isometries (see \cite[Proposition~2.11]{KP}). In a $k$-graph
which admits sources, however, not every finite path is
contained in an infinite path of the form defined in
\cite{KP}, and hence the proof of
\cite[Proposition~2.11]{KP} does not carry over. To allow
for sources, we replace $\Lambda^\infty$ with a
boundary-path space $\uptoLambda(\infty)$; for locally
convex $k$-graphs, we can achieve this construction using
the $k$-graphs $\Omega_{k,m}$ of Example~\ref{k-graphs}(ii).

\begin{dfn} \label{def:uptoLambda(infty)}
Let $\Lambda$ be a locally convex $k$-graph. A {\it boundary
path} in $\Lambda$ is a graph mor\-phism $x : \Omega_{k,m}
\to \Lambda$ for some $m\in(\NN\cup\{\infty\})^k$ such that,
whenever $v\in\Obj(\Omega_{k,m})$ satisfies
$(\Omega_{k,m})^{\le e_i}(v) = \{v\},$ we also have
$\uptoLambda(e_i)(x(v)) =\{x(v)\}.$ We denote the collection
of all boundary paths in $\Lambda$ by $\uptoLambda(\infty).$
The range map of $\Lambda$ extends naturally to
$\uptoLambda(\infty)$ via $r(x) := x(0).$ For $v \in
\Lambda^0,$ we write $\uptoLambda(\infty)(v)$ for $\{x \in
\uptoLambda(\infty) : r(x) = v\}.$ We define a degree map
$d_\infty:\uptoLambda(\infty)\to (\NN\cup\{\infty\})^k$ by
setting $d_\infty(x):=m$ when $x:\Omega_{k,m}\to\Lambda.$
\end{dfn}

As with the infinite paths of \cite{KP}, a boundary path $x$
is completely determined by the set of paths $\{x(0,p) : p
\le d_\infty(x)\}.$ In fact, when the $k$-graph has no
sources, $\uptoLambda(\infty)$ is exactly the infinite path
space from \cite{KP}. If a $k$-graph $\Lambda$ is locally
convex then for any vertex $v$ of $\Lambda,$ the set
$\uptoLambda(\infty)(v)$ is nonempty: even if $v$ emits no
paths of nonzero degree, we have $\uptoLambda(\infty)(v) =
\{v\} \not= \emptyset.$

\begin{theorem} \label{spadesuit theorem}
Let $(\Lambda,d)$ be a row-finite $k$-graph. Then there is a
Cuntz-Krieger $\Lambda$-family
$\{S_\lambda:\lambda\in\Lambda\}$ with each $S_\lambda$
nonzero if and only if $\Lambda$ is locally convex.
\end{theorem}

\begin{proof}
First suppose that $\Lambda$ is not locally convex. Then
there exists a vertex $v\in\Lambda^0$ and
$\mu\in\Lambda^{e_i}(v)$ for some $i\in\{1,\dots,k\}$ such
that $\Lambda^{e_j}(v)\neq\emptyset$ and
$\Lambda^{e_j}(s(\mu))=\emptyset$ for some $j\neq i$.
Considering the partial isometry $s_\mu\in C^*(\Lambda)$, we
have
\[
s_\mu = s_vs_\mu = \sum_{\nu\in\Lambda^{e_j}(v)}s_\nu
s^*_\nu s_\mu
= \sum_{\nu\in\Lambda^{e_j}(v)}s_\nu
\sum_{\substack{\nu\alpha=\mu\beta \\
d(\mu\beta)=e_i+e_j}} s_\alpha s^*_\beta,
\]
but since $\Lambda^{e_j}(s(\mu))=\emptyset$, no such $\beta$
exists. Thus $s_\mu=0$, and so by the universal property of
$C^*(\Lambda)$ any Cuntz-Krieger $\Lambda$-family
$\{S_\lambda:\lambda\in\Lambda\}$ must have $S_\mu=0$.

Now suppose that $\Lambda$ is locally convex. Let
$\Hh:=\ell^2(\uptoLambda(\infty))$, and for each
$\lambda\in\Lambda$ define $S_\lambda\in B(\Hh)$ by
\begin{equation}\label{S_lambda eqn}
S_\lambda u_x :=
\begin{cases}
u_{\lambda x} \quad\text{if $s(\lambda)=r(x)$} \\
0 \quad\text{otherwise},
\end{cases}
\end{equation}
where $\{u_x:x\in\uptoLambda(\infty)\}$ is the usual basis
for $\Hh$. Each $S_\lambda \not= 0$ because
$\uptoLambda(\infty)(s(\lambda))\neq\emptyset$.
Cuntz-Krieger relations (1)-(3) follow directly from the
definition of the operators $S_\lambda$; it remains to show
that relation~(4) is fulfilled. Since $\Lambda$ is locally
convex, by Proposition~\ref{equiv rel} it suffices to show
that for each $v\in\Lambda^0$ and $i\in\{1,\dots,k\}$,
$S_v=\displaystyle{\sum_{\lambda\in\uptoLambda(e_i)(v)}S_\lambda S^*_\lambda}$. If
$\uptoLambda(e_i)(v)=\{v\}$, then the relation is trivially true, so suppose
$\uptoLambda(e_i)(v)\neq\{v\}$. For
$x\in\uptoLambda(\infty)$ we have
\begin{align*}
\sum_{\lambda\in\uptoLambda(e_i)(v)}S_\lambda S^*_\lambda
u_x &=
\sum_{\lambda\in\uptoLambda(e_i)(v)}S_\lambda
(\delta_{\lambda,x(0,e_i)}u_{x(e_i, \infty)}) \\
&=\sum_{\lambda\in\uptoLambda(e_i)(v)}\delta_{\lambda,x(0,e_i)}u_x \\
&=\begin{cases}u_x \quad\text{if there exists
$\lambda\in\Lambda^{e_i}(v)$
such that $\lambda=x(0,e_i)$} \\
0 \quad\text{otherwise}.
\end{cases}
\end{align*}
Taking $\lambda = v$ in \eqref{S_lambda eqn}, we can see
that it suffices to show that $r(x)=v$ if and only if there
exists $\lambda\in\Lambda^{e_i}(v)$ such that
$\lambda=x(0,e_i)$. If $\lambda=x(0,e_i)$ for some
$\lambda\in\Lambda^{e_i}(v)$, then $r(x)=r(\lambda)=v$. If
$r(x)=v$, then $x(0,e_i)\in\Lambda^{e_i}(v)$ because
$\uptoLambda(e_i)(v)\neq\{v\}$ and $x$ is a boundary path.
\end{proof}

\section{The uniqueness theorems}

\subsection{The gauge-invariant uniqueness
theorem}\label{grunt}  Our gauge-invariant uniqueness
theorem extends \cite[Theorem~3.4]{KP} to row-finite
$k$-graphs with sources.

Let $(\Lambda, d)$ be a row-finite $k$-graph. For
$z\in\TT^k$ and $n \in \ZZ^k$, let $z^n:=z_1^{n_1} \cdots
z_k^{n_k}.$ Then
$\{z^{d(\lambda)}s_\lambda:\lambda\in\Lambda\}$ is a
Cuntz-Krieger $\Lambda$-family which generates
$C^*(\Lambda)$, and the universal property of $C^*(\Lambda)$
gives a homomorphism $\gamma_z:C^*(\Lambda)\to C^*(\Lambda)$
such that $\gamma_z(s_\lambda)=z^{d(\lambda)}s_\lambda$ for
$\lambda\in\Lambda$; $\gamma_{\bar z}$ is an inverse for
$\gamma_z,$ so $\gamma_z$ is an automorphism. An $\epsilon/3$-argument shows that $\gamma$ is
a strongly continuous action of $\TT^k$ on $C^*(\Lambda)$, which is called the
{\it gauge action}.

\begin{theorem}[The Gauge-Invariant Uniqueness Theorem]
\label{giut}
Let $(\Lambda,d)$ be a locally convex row-finite $k$-graph,
let $\{t_\lambda:\lambda\in\Lambda\}$ be a Cuntz-Krieger
$\Lambda$-family, and let $\pi$ be the representation of
$C^*(\Lambda)$ such that $\pi(s_\lambda)=t_\lambda$. If each
$t_v$ is nonzero and there is a strongly continuous action
$\beta:\TT^k\to \Aut(C^*(\{t_\lambda:\lambda\in\Lambda\}))$
such that $\beta_z\circ\pi=\pi\circ\gamma_z$ for
$z\in\TT^k$, then $\pi$ is faithful.
\end{theorem}

\begin{rmk}
Strictly speaking, it is not necessary to assume that
$\Lambda$ is locally convex; if there is a Cuntz-Krieger
$\Lambda$-family $\{t_\lambda : \lambda \in \Lambda\}$ with
each $t_v$ nonzero, then Theorem~\ref{spadesuit theorem}
implies that $\Lambda$ is locally convex.
\end{rmk}

The first part of the proof, the analysis of the core
$C^*(\Lambda)^\gamma$, is the same for both uniqueness
theorems. Using our$\uptoLambda(p)$, and Lemmas~\ref{lambda
inclusion} and \ref{orth range proj}, we can follow the
argument of \cite[\S3]{KP}. We consider the map $\Phi:
C^*(\Lambda) \to C^*(\Lambda)$ defined by
\[
\Phi(a) := \int_{\TT^k}\gamma_z(a)\,dz,
\]
which is faithful on positive elements, and has range the
fixed point algebra $C^*(\Lambda)^\gamma.$ To identify the
structure of $C^*(\Lambda)^\gamma$, we let $v \in
\Lambda^0$, $q \in \NN^k$, and define
\[
\Ff_q(v):=\clsp\{s_\lambda s^*_\mu:
\lambda,\mu\in\uptoLambda(q),
d(\lambda)=d(\mu),s(\lambda)=s(\mu)=v\}.
\]
It follows from Lemma~\ref{orth range proj} that $\Ff_q(v)$
is the direct sum of the subalgebras
\[
\Ff_{q,p}(v):=\clsp\{s_\lambda s^*_\mu:
\lambda,\mu\in\uptoLambda(q),
d(\lambda)=d(\mu)=p,s(\lambda)=s(\mu)=v\},
\]
that
\begin{equation}\label{isowithK}
\Ff_{q,p}(v)=\Kk(\ell^2(\{\lambda\in\uptoLambda(q):d(\lambda
)=p \text{ and } s(\lambda)=v\})),
\end{equation}
and that for fixed $q$, the $\Ff_q(v)$'s are mutually
orthogonal. Since the elements $s_\lambda s^*_\mu$ span
$C^*(\Lambda)$, and since $\gamma_z(s_\lambda s^*_\mu) =
z^{d(\lambda)-d(\mu)} s_\lambda s^*_\mu$, the algebras
$\Ff_q := \oplus_{v \in \Lambda^0} \Ff_q(v)$ span
$C^*(\Lambda)^\gamma$. When $ p \le q$, we have $\Ff_p
\subset \Ff_q$ by Lemma~\ref{lambda inclusion}, so
$C^*(\Lambda)^\gamma$ is the direct limit $\overline{\cup_{q
\in \NN^k} \Ff_q}$ of the algebras $\Ff_q$. In particular,
$C^*(\Lambda)^\gamma$ is AF.

\begin{proof}[Proof of Theorem~\ref{giut}]
Because each $t_v$ is nonzero, and $t_\lambda$ has initial
projection $t_{s(\lambda)}$, each $t_\lambda$ is nonzero,
and hence the representation $\pi$ is nonzero on each
$\Ff_{q,p}(v)$. It follows from (\ref{isowithK}) that $\pi$
is faithful on $\Ff_{p,q}(v)$, hence on $\Ff_q(v)$ and on
$\Ff_q$. Since $C^*(\Lambda)^\gamma=\varinjlim\Ff_q$, it
follows that $\pi$ is faithful on $C^*(\Lambda)^\gamma$ (see
\cite[Lemma~1.3]{ALNR}, for example). We can now use the
argument of \cite[page~11]{KP}.
\end{proof}

\subsection{The Cuntz-Krieger uniqueness theorem} Our
Cuntz-Krieger uni\-queness theorem extends
\cite[Theorem~4.6]{KP} to row-finite $k$-graphs with sources
(see Remark~\ref{aperiodicity} below).

\begin{theorem} \label{CKUThm}
Let $(\Lambda,d)$ be a locally convex row-finite $k$-graph,
and suppose that:
\begin{equation}\label{B}
\text{for each $v \in \Lambda^0$, there exists $x \in
\uptoLambda(\infty)(v)$ such that $\alpha \not= \beta$
implies $\alpha x \not= \beta x$.}\tag{B}
\end{equation}
Let $\{t_\lambda:\lambda\in\Lambda\}$ be a Cuntz-Krieger
$\Lambda$-family, and let $\pi$ be the representation of
$C^*(\Lambda)$ such that $\pi(s_\lambda) = t_\lambda$. If
each $t_v$ is non-zero, then $\pi$ is faithful.
\end{theorem}

\begin{proof}
For $\lambda \in \Lambda \cup \uptoLambda(\infty)$ and $p
\in \NN^k$, we write $\lambda(0,p)$ for the unique path such
that $\lambda=\lambda(0,p)\lambda'$ with $d(\lambda(0,p))=
p\wedge d(\lambda)$.

We know from the analysis in \S\ref{grunt} that $\pi$ is
faithful on the fixed-point algebra $C^*(\Lambda)^\gamma$,
and hence the standard argument will work once we know that
\begin{equation} \label{CKU norm inequality}
\|\pi(\Phi(a))\|\le \|\pi(a)\| \text{ for $a\in
C^*(\Lambda)$}.
\end{equation}
Recalling that $\lsp\{s_\lambda s^*_\mu:
s(\lambda)=s(\mu)\}$ is dense in $C^*(\Lambda)$, we consider
arbitrary $a=\sum_{(\lambda,\mu)\in
F}\zeta_{(\lambda,\mu)}s_\lambda s^*_\mu$ where $F$ is
finite and $\zeta_{(\lambda,\mu)}\in\CC$. Let $l$ be the
least upper bound of $\{d(\lambda) \vee d(\mu) :
(\lambda,\mu) \in F\}.$ Then
\[
\Phi(a)=\sum_{\{(\lambda,\mu)\in
F\,:\,d(\lambda)=d(\mu)\}}\zeta_{(\lambda,\mu)}s_\lambda
s^*_\mu \in\Ff_l.
\]
For $(\lambda,\mu)\in F$, we define
\[
F_{(\lambda,\mu)}=\{(\lambda\nu,\mu\nu):\nu\in\uptoLambda({l
-d(\lambda)})(s(\lambda))\},
\]
and $E=\cup_{(\lambda,\mu)\in F}F_{(\lambda,\mu)}$. For $\nu
\in\uptoLambda({l-d(\lambda)})(s(\lambda))$, we define
\[
\xi_{(\lambda\nu,\mu\nu)}=\sum_{\{(\lambda^\prime,\mu^\prime
)\in F\,:\,(\lambda\nu,\mu\nu)=(\lambda^\prime\nu^\prime,
\mu^\prime\nu^\prime) \text{ for some }
\nu^\prime\in\Lambda\}}\zeta_{(\lambda^\prime,\mu^\prime)},
\]
and using Cuntz-Krieger relation (4) we then have
\[
a=\sum_{(\alpha,\beta)\in E}\xi_{(\alpha,\beta)}s_\alpha
s^*_\beta;
\]
the point is that now $\alpha\in\uptoLambda(l)$ for all
$(\alpha,\beta)\in E$.

Since $\Ff_l$ decomposes as a direct sum
$\oplus_{v\in\Lambda^0}\Ff_l(v)$, so does its image under
$\pi$, and there is a vertex $v\in\Lambda^0$ such that
\begin{equation}\label{norm attained}
\|\pi(\Phi(a))\|=\Big\|\sum_{\{(\alpha,\beta)\in
E\,:\,d(\alpha)=d(\beta)
\text{ and } s(\alpha) = v\}} \xi_{(\alpha,\beta)} t_\alpha
t^*_\beta\Big\|.
\end{equation}
Choose a boundary path $x\in\uptoLambda(\infty)(v)$ such
that $\alpha x\neq \beta x$ for all $\alpha \not=
\beta\in\Lambda$; then for each $\alpha \not=
\beta\in\Lambda$, there exists $M_{\alpha,\beta} \ge
d(\alpha) \vee d(\beta)$ such that $(\alpha x)(0,m)\neq
(\beta x)(0,m)$ whenever $m \ge M_{\alpha,\beta}$. Let
\[
T = \{\tau \in\uptoLambda(l) : \tau=\alpha\text{ or
}\tau=\beta\text{ for some }(\alpha,\beta) \in E, s(\tau) =
v\}.
\]
Let $M$ be the least upper bound of $\{M_{\tau, \beta} :
\tau \in T, (\alpha,\beta)\in E\text{ for some }\alpha\}$.
In particular,
\begin{equation} \label{aperiodic path}
(\beta x)(0,M) \neq (\tau x)(0,M)
\end{equation}
when $\beta$ is the second coordinate of an element of $E$,
$\tau \in T$, and $\beta \not= \tau$. Write $x_M$ for
$x(0,M)$.

For each $n \le l,$ we define
\[
Q_n := \sum_{\{\tau\in T\,:\,d(\tau)=n\}}t_{\tau
x_M}t^*_{\tau x_M}\Bigg.
\]
Now we define $Q : C^*(\{t_\lambda : \lambda \in
\Lambda\})\to C^*(\{t_\lambda : \lambda \in \Lambda\})$ by
\[
Q(b):= \sum_{n \le l} Q_n b Q_n.
\]
Since the $Q_n$ are mutually orthogonal projections, we have
\[
\|Q(b)\| = \Big\|\sum_{n \le l} Q_n b Q_n\Big\| \le \|b\|
\text{ for $b \in C^*(\{t_\lambda\})$}.
\]
We aim to show that $\|Q(\pi(\Phi(a)))\| =
\|\pi(\Phi(a))\|,$ and that $Q(\pi(a))=Q(\pi(\Phi(a)))$;
this will give us
\begin{equation}\label{todo}
\|\pi(\Phi(a))\| = \|Q(\pi(\Phi(a)))\| = \|Q(\pi(a))\| \le
\|\pi(a)\|,
\end{equation}
and the proof will be complete.

Write $M_T$ for the matrix algebra spanned by $\{s_\lambda
s^*_\mu : \lambda,\mu \in T, d(\lambda) = d(\mu)\}$. Notice
that $M_T \subset \Ff_l(v)$. For $s_\lambda s^*_\mu \in M_T$
we have $\lambda, \mu \in \uptoLambda(l),$ so for $\tau \in
T,$ $t^*_\tau t_\lambda = 0$ unless $\tau = \lambda,$ and
$t^*_\mu t_\tau = 0$ unless $\tau = \mu,$ and hence
\begin{align*}
Q(\pi(s_\lambda s^*_\mu))
&= \sum_{n\le l}\Bigg(\sum_{\{\tau\in
T\,:\,d(\tau)=n\}}t_{\tau x_M}t^*_{\tau x_M}\Bigg)
t_\lambda t^*_\mu \Bigg(\sum_{\{\tau'\in
T\,:\,d(\tau')=n\}}t_{\tau' x_M}t^*_{\tau' x_M}\Bigg) \\
&= t_{\lambda x_M} t^*_{x_M} t^*_\lambda t_\lambda t^*_\mu
t_\mu t_{x_M} t^*_{\mu x_M} \\
&= t_{\lambda x_M} t^*_{\mu x_M} \\
&\neq 0.
\end{align*}
Using Lemma~\ref{orth range proj}, it follows that $\{
Q(\pi(s_\lambda s^*_\mu)) : s_\lambda s^*_\mu\in M_T\}$ is a
family of nonzero matrix units, and from this we deduce that
the map $b \mapsto Q(\pi(b))$ is a faithful representation
of $M_T.$ Since both $\pi$ and $Q \circ \pi$ are faithful on
$M_T$ and since
\[
\sum_{\{(\alpha,\beta)\in E\,:\,d(\alpha)=d(\beta) \text{
and }
s(\alpha) = v\}} \xi_{(\alpha,\beta)} t_\alpha t^*_\beta \in
M_T,
\]
it follows from \eqref{norm attained} that $\|\pi(\Phi(a))\|
= \|Q(\pi(\Phi(a)))\|$.

To establish that $Q(\pi(a))=Q(\pi(\Phi(a)))$, we show that
$Q(t_\alpha t^*_\beta)=0$ whenever $(\alpha,\beta)\in E$ and
$d(\alpha) \neq d(\beta);$ this shows that $Q$ kills those
terms of $\pi(a)$ which are the images under $\pi$ of terms
of $a$ killed by $\Phi$. Notice that if $(\alpha, \beta) \in
E$ then $\alpha \in \uptoLambda(l),$ so for $\tau \in T,$
$t^*_\tau t_\alpha = 0$ unless $\tau = \alpha.$ Hence, for
$(\alpha,\beta)\in E$ with $d(\alpha)\neq d(\beta)$, we have
\begin{align*}
Q(t_\alpha t^*_\beta)
&= \sum_{n\le l}\Bigg(\sum_{\{\tau\in
T\,:\,d(\tau)=n\}}t_{\tau x_M}t^*_{\tau x_M}\Bigg)
t_\alpha t^*_\beta \Bigg(\sum_{\{\tau'\in
T\,:\,d(\tau')=n\}}t_{\tau' x_M}t^*_{\tau' x_M}\Bigg) \\
&= \sum_{\{\tau' \in T\,:\,d(\tau')=d(\alpha)\}} t_{\alpha
x_M} t^*_{\beta x_M} t_{\tau' x_M} t^*_{\tau' x_M} \\
&= \sum_{\{\tau' \in T\,:\,d(\tau')=d(\alpha)\}} t_{\alpha
x_M}
        \Bigg(\sum_{\substack{\beta x_M \eta = \tau' x_M
\zeta \\ d(\beta x_M\eta) = d(\beta x_M) \vee d(\tau' x_M)}}
                  t_\eta t^*_\zeta\Bigg)
        t^*_{\tau' x_M},
\end{align*}
which is nonzero if and only if there exist $\eta,\zeta \in
\Lambda$ such that
\begin{equation}\label{paths equal}
(\beta x_M \eta)(0,M) = (\tau' x_M \zeta)(0,M).
\end{equation}
But $(\beta x_M \eta)(0,M) = (\beta x_M)(0,M)$ : if not,
then there exists an $i$ such that $d(\beta x_M)_i < M_i$
and $d(\eta)_i > 0.$ But
\[
d(\beta x_M)_i < M_i \implies d(x_M)_i < M_i \implies
\Lambda^{e_i}(s(x_M)) = \emptyset
\]
since $x$ is a boundary path. Likewise, $(\tau' x_M
\zeta)(0,M) = (\tau' x_M)(0,M)$, and so \eqref{paths equal}
is equivalent to $(\beta x_M)(0,M) = (\tau' x_M)(0,M)$,
which is impossible by \eqref{aperiodic path}. This proves
(\ref{todo}), and the result follows.
\end{proof}

\begin{rmk}\label{aperiodicity}
The condition \eqref{B} in Theorem~\ref{CKUThm} is automatic
if $\Lambda$ has no sources and satisfies the aperiodicity
condition (A) of \cite[Definition~4.3]{KP}. To see this, let
$\sigma$ be the shift map on $\Lambda^\infty$ defined as in
\cite{KP} by $\sigma^p(x) = x(p,\infty)$. Suppose that
$\Lambda$ has no sources and \eqref{B} does not hold. Then
there is a vertex $v \in \Lambda^0$ such that for each $x
\in \Lambda^\infty(v)$, there exist $\alpha_x \not= \beta_x$
such that $\alpha_x x = \beta_x x$. Then for $x \in
\Lambda^\infty(v)$,
\[
\begin{split}
\sigma^{d(\alpha_x) \vee d(\beta_x) - d(\alpha_x)}(x) &=
\sigma^{d(\alpha_x) \vee d(\beta_x)}(\alpha_x x) \\
&= \sigma^{d(\alpha_x) \vee d(\beta_x)}(\beta_x x) =
\sigma^{d(\alpha_x) \vee d(\beta_x) - d(\beta_x)}(x).
\end{split}
\]
Hence condition (A) of \cite[Definition~4.3]{KP} does not
hold at $v$.

Thus Theorem~\ref{CKUThm} is formally stronger than
\cite[Theorem~4.6]{KP} even when $\Lambda$ has no sources.
We have been unable to decide whether it is equivalent: we
do not know whether in graphs without sources, \eqref{B}
implies (A). For 1-graphs with no sources, we can prove that
\eqref{B} implies (A) because it is easy to construct
aperiodic paths (see the proof of \cite[Lemma~3.4]{KPR}). In
higher-rank graphs, it is hard to produce aperiodic paths,
and we suspect that in practice \eqref{B} might be easier to
check than (A).
\end{rmk}

\section{The ideal structure}
Let $(\Lambda,d)$ be a locally convex row-finite $k$-graph.
Define a relation on $\Lambda^0$ by setting $v\ge w$ if
there is a path $\lambda\in\Lambda$ with $r(\lambda)=v$ and
$s(\lambda)=w$. A subset $H$ of $\Lambda^0$ is {\it
hereditary} if $v\ge w$ and $v\in H$ imply $w\in H$; $H$ is
{\it saturated} if for $v\in\Lambda^0$,
\[
\{s(\lambda): \lambda\in\uptoLambda({e_i})(v)\}\subset H
\text{~for some~} i\in\{1,\dots,k\} \implies v\in H.
\]
The {\it saturation} of a set $H$ is the smallest saturated
subset ${\overline H}$ of $\Lambda^0$ containing $H$.

\begin{lemma} \label{sat hered}
Suppose $\Lambda$ is a locally convex row-finite $k$-graph,
and $H$ is a hereditary subset of $\Lambda^0$. Then the
saturation $\overline{H}$ is hereditary.
\end{lemma}
\begin{proof}
We use an inductive construction of $\overline{H}$ like that
used by Szyma\`nski for 1-graphs in \cite{S01}. For $F
\subset \Lambda^0$, we define
\[
\Sigma(F) := \bigcup^k_{i=1}\{v \in \Lambda^0 : s(\lambda)
\in F\text{ for all }\lambda \in \uptoLambda(e_i)(v)\},
\]
and write $\Sigma^n(F)$ for the set obtained by repeating
the process $n$ times. Notice that if $F$ is hereditary,
then $F \subset \Sigma(F)$. We will show that
$\bigcup^\infty_{n=1} \Sigma^n(H)$ is hereditary and equal
to $\overline{H}$.

We begin by showing that if $F$ is hereditary, then
$\Sigma(F)$ is hereditary. To see this, suppose that $v \in
\Sigma(F)$ and that $v \ge w.$ Then there exists $\lambda
\in \Lambda^0$ such that $r(\lambda) = v$ and $s(\lambda) =
w.$ If $d(\lambda) = 0$, then $w = v \in F$, so suppose
$d(\lambda)_j > 0$, and factor $\lambda = \lambda^\prime
\lambda^{\prime\prime}$ where $d(\lambda^\prime) = e_j$. We
claim that $s(\lambda^\prime) \in \Sigma(F)$. To see this,
choose $i$ such that $\{s(\mu) : \mu \in
\uptoLambda(e_i)(v)\} \subset F.$ If $\uptoLambda(e_i)(v) =
\{v\}$ or if $i = j$, then $s(\lambda^\prime) \in F \subset
\Sigma(F)$. So suppose that $\uptoLambda(e_i)(v) \not=
\{v\}$ and $i \not= j$. Since $\Lambda$ is locally convex,
$\Lambda^{e_i}(s(\lambda^\prime)) \not= \emptyset$, so it
suffices to show that $\nu \in
\Lambda^{e_i}(s(\lambda^\prime))$ implies $s(\nu) \in F$.
Let $\nu \in \Lambda^{e_i}(s(\lambda^\prime))$. Then
$\lambda^\prime\nu = \mu\nu^\prime$ for some $\mu \in
\uptoLambda(e_i)(v)$. But now $s(\mu) \in F$ and hence
$s(\nu^\prime) = s(\nu) \in F$ because $F$ is hereditary.
Thus $s(\lambda^\prime) \in \Sigma(F)$ as claimed. By
induction on the length of $\lambda$, it follows that $w \in
\Sigma(F)$, and hence $\Sigma(F)$ is hereditary.

We now know that $\Sigma^n(H) \subset \Sigma^{n+1}(H)$ for
all $n$, and that $\Sigma^n(H)$ is hereditary for all $n$;
thus $\bigcup^\infty_{n=1} \Sigma^n(H)$ is also hereditary.
It remains to show that $\overline{H} = \bigcup^\infty_{n=1}
\Sigma^n(H)$. Because applying $\Sigma$ can never take us
outside of a saturated set, we have $\bigcup\Sigma^n(H)
\subset \overline{H}$, so it is enough to show that
$\bigcup^\infty_{n=1} \Sigma^n(H)$ is saturated. To see
this, suppose that $v \in \Lambda^0$ and $\{s(\lambda) :
\lambda \in \uptoLambda(e_i)(v)\} \subset
\bigcup^\infty_{n=1} \Sigma^n(H)$. Then, since $\Lambda$ is
row finite, we have $\{s(\lambda) : \lambda \in
\uptoLambda(e_i)(v)\} \subset \Sigma^N(H)$ for some $N$ and
it follows that $v \in \Sigma^{N+1}(H)$. Thus
$\bigcup^\infty_{n=1} \Sigma^n(H)$ is saturated.
\end{proof}

\begin{theorem} \label{ideals}
Let $(\Lambda, d)$ be a locally convex row-finite $k$-graph.
For each subset $H$ of $\Lambda^0$, let $I_H$ be the closed
ideal in $C^*(\Lambda)$ generated by $\{s_v : v \in H\}$.
\begin{itemize}
\item[(a)] The map $H \mapsto I_H$ is an isomorphism of the
lattice of saturated hereditary subsets of $\Lambda^0$ onto
the lattice of closed gauge-invariant ideals of
$C^*(\Lambda)$.
\item[(b)] Suppose $H$ is saturated and hereditary. Then
\[
\Gamma(\Lambda \setminus H) := (\Lambda^0 \setminus H,
\{\lambda \in \Lambda : s(\lambda) \not \in H\}, r, s)
\]
is a locally convex row-finite $k$-graph, and $C^*(\Lambda)
/ I_H$ is canonically isomorphic to $C^*(\Gamma(\Lambda
\setminus H))$.
\item[(c)] If $H$ is any hereditary subset of $\Lambda^0$,
then
\[
\Lambda(H) := (H, \{\lambda \in \Lambda : r(\lambda) \in
H\}, r, s)
\]
is a locally convex row-finite $k$-graph, $C^*(\Lambda(H))$
is canonically isomorphic to the subalgebra $C^*(\{s_\lambda
: r(\lambda) \in H\})$ of $C^*(\Lambda)$, and this
subalgebra is a full corner in $I_H$.
\end{itemize}
\end{theorem}
\begin{proof}
The proof of Theorem~\ref{ideals} is the same as the proof
of \cite[Theorem~4.1]{BPRS} once we establish that
$\Gamma(\Lambda \setminus H)$ and $\Lambda(H)$ from parts
(b) and (c) are locally convex row-finite $k$-graphs. This
is easy to check for $\Lambda(H)$, and the row-finiteness of
$\Gamma(\Lambda \setminus H)$ follows from that of
$\Lambda$. We need to check that $\Gamma(\Lambda \setminus
H)$ is a $k$-graph and is locally convex. For convenience,
write $\Gamma = \Gamma(\Lambda \setminus H)$.

To show that the factorisation property holds for
$(\Gamma,d)$, take $\lambda\in\Gamma$ and suppose
$d(\lambda)=p+q$. We know there exist unique
$\mu,\nu\in\Lambda$ such that $\lambda=\mu\nu$, $d(\mu)=p$
and $d(\nu)=q$. Certainly $s(\nu)=s(\lambda)\notin H$, and
if $s(\mu)\in H$, then $s(\nu)\in H$, a contradiction. Hence
$\mu,\nu\in\Gamma$, and $\Gamma$ is a sub-$k$-graph.

Now we show that $\Gamma$ is locally convex. Consider an
arbitrary vertex $v\in\Gamma^0$ which has
$\lambda\in\Gamma^{e_i}(v)$ and $\mu\in\Gamma^{e_j}(v)$ for
some $i\neq j$. We know that $s(\lambda),s(\mu)\notin H$,
and since $\Lambda$ is locally convex, we also know there
exist $\alpha\in\Lambda^{e_i}(s(\mu))$ and
$\beta\in\Lambda^{e_j}(s(\lambda))$. If
$\Gamma^{e_i}(s(\mu))=\emptyset$, then
$\{s(\alpha):\alpha\in\uptoLambda({e_i})(s(\mu))\}\subset
H$, and similarly, if $\Gamma^{e_j}(s(\mu))=\emptyset$, then
$\{s(\beta):\beta\in\uptoLambda({e_j})(s(\lambda))\}\subset
H$; in either case saturatedness implies that $s(\mu)\in H$
or $s(\lambda)\in H$, a contradiction. Hence
$\Gamma^{e_i}(s(\mu))\neq\emptyset$ and
$\Gamma^{e_j}(s(\lambda))\neq\emptyset$, so $\Gamma$ is
locally convex.
\end{proof}

The proof of the next theorem is the same as the first two
paragraphs of the proof of \cite[Theorem~4.1]{BPRS} except
that in the first paragraph, we apply the Cuntz-Krieger
uniqueness theorem rather than the gauge-invariant
uniqueness theorem to show that $H \mapsto I_H$ is onto.

\begin{theorem}
Let $(\Lambda,d)$ be a locally convex row-finite $k$-graph
such that for every saturated hereditary subset $H$ of
$\Lambda^0$, $\Gamma(\Lambda \setminus H)$ satisfies
condition \eqref{B} of Theorem~\ref{CKUThm}. Then $H\mapsto
I_H$ is an isomorphism of the lattice of saturated
hereditary subsets of $\Lambda^0$ onto the lattice of closed
ideals of $C^*(\Lambda)$.
\end{theorem}


\begin{thebibliography}{00}
\bibitem{ALNR} S. Adji, M. Laca, M. Nilsen, and I. Raeburn,
\emph{Crossed products by semigroups of endomorphisms and
the Toeplitz algebras of ordered groups}, Proc. Amer. Math.
Soc. {\bf 122} (1994), 1133-1141.

\bibitem{B} B. Blackadar, \emph{Shape theory for
$C^*$-algebras}, Math. Scand. {\bf 56} (1985), 249--275.

\bibitem{BPRS} T. Bates, D. Pask, I. Raeburn, and W.
Szyma\'nski, \emph{The $C^*$-algebras of row--finite
graphs}, New York J. Math. {\bf 6} (2000), 307--324.

\bibitem{C} J. Cuntz, \emph{A class of $C^*$-algebras and
topological Markov chains II: reducible chains and the
Ext-functor for $C^*$-algebras}, Invent. Math. {\bf 63}
(1981), 25--40.

\bibitem{CK} J. Cuntz and W. Krieger, \emph{A class of
$C^*$-algebras and topological Markov chains}, Invent. Math.
{\bf 56} (1980), 251--268.

\bibitem{FS} N. J. Fowler and A. Sims, \emph{Product systems
over right-angled Artin semigroups}, Trans. Amer. Math.
Soc., to appear.

\bibitem{KP} A. Kumjian and D. Pask, \emph{Higher rank graph
$C^*$-algebras}, New York J. Math. {\bf 6} (2000), 1--20.

\bibitem{KPR} A. Kumjian, D. Pask, and I. Raeburn,
\emph{Cuntz-Krieger algebras of directed graphs}, Pacific J.
Math. {\bf 184} (1998), 161--174.

\bibitem{KPRR} A. Kumjian, D. Pask, I. Raeburn, and J.
Renault, \emph{Graphs, groupoids and Cuntz-Krieger
algebras}, J. Funct. Anal. {\bf 144} (1997), 505--541.

\bibitem{RW} I. Raeburn and D. P. Williams, \emph{Morita
equivalence and continuous-trace $C^*$-algebras}, Math.
Surveys \& Monographs, vol. 60, Amer. Math. Soc. Providence,
1998.

\bibitem{RS1} G. Robertson and T. Steger,
\emph{$C^*$-algebras arising from group actions on the
boundary of a triangle building}, Proc. London Math. Soc.
{\bf 72} (1996), 613--637.

\bibitem{RS2} G. Robertson and T. Steger, \emph{Affine
buildings, tiling systems and higher rank Cuntz-Krieger
algebras}, J. reine angew. Math. {\bf 513} (1999), 115--144.

\bibitem{RS3} G. Robertson and T. Steger, \emph{Asymptotic
K-theory for groups acting on $\tilde A_2$ buildings},
Canad. J. Math., to appear.

\bibitem{S01} W. Szyma\'nski, \emph{Simplicity of
Cuntz-Krieger algebras of infinite matrices}, Pacific J.
Math. {\bf 199} (2001), 249--256.
\end{thebibliography}
\end{document}